\documentclass[reqno, a4paper, 12pt]{amsart}
\usepackage[utf8]{inputenc}
\usepackage[english]{babel}
\usepackage{csquotes}

\usepackage{amsfonts, amsmath, amssymb, amsthm,caption}
\usepackage{aliascnt} % for cleveref-aware theorem environments sharing counters
\usepackage{mathtools}
\usepackage{mathrsfs}
\usepackage{bbm}
\usepackage{breqn}
\usepackage{bm}
\usepackage{cancel}
\usepackage{comment}

\usepackage{fullpage}

\usepackage[url=false, isbn=false, eprint=false, maxbibnames=99, doi=false,  maxcitenames=99]{biblatex}
\usepackage{ragged2e}
\usepackage{enumerate}
\usepackage{enumitem}
\usepackage{multicol}
\usepackage{tasks}
\usepackage{hyperref}
\usepackage[left=1in,right=1in,top=1in,bottom=1in]{geometry}

\usepackage{graphicx} % Required for inserting images
\usepackage{tikz,mathtools,mathrsfs}
\usepackage{xcolor}
\newtheorem{theorem}{Theorem}[section]
\usepackage{cleveref}

\newaliascnt{lemma}{theorem}
\newtheorem{lemma}[lemma]{Lemma}
\aliascntresetthe{lemma}

\newtheorem{example}[theorem]{Example}

\newtheorem{question}[theorem]{Question}

\crefname{lemma}{Lemma}{Lemmas}
\Crefname{lemma}{Lemma}{Lemmas}

\newcommand{\pmcggg}{{{\mathcal{M}_p\left(\mathcal{G}_{3\times 2n}\right)}}}

\DeclareMathOperator{\st}{st}
\DeclareMathOperator{\del}{del}
\DeclareMathOperator{\lk}{lk}

\title{Perfect Matching Complexes of $3\times 2n$ Grid Graph}

\author{ Himanshu Chandrakar }

\address{Chennai Mathematical Institute, Siruseri, Tamil Nadu, India.}

\email{himanshuc@cmi.ac.in} 

\thanks{Chandrakar is partially supported by a grant from the Infosys Foundation.}
\author{Trung Chau}

\address{University of Manitoba, Winnipeg, Manitoba, Canada.}

\email{chauchitrung1996@gmail.com}

\thanks{Chau is partially supported in part by the Pacific Institute for the Mathematical Sciences and a grant from the Simons Foundation International [SFI-MPS-T-Institutes-00020822, OY]}

\author{Shuchita Goyal}

\address{Birla Institute of Technology and Science Pilani, Pilani Campus, Rajasthan, India} 

\email{shuchita.goyal@pilani.bits-pilani.ac.in}
\thanks{Goyal is partially supported by the NFSG grant N5/25/1045}

\author{Rahul Karki}

\address{UPES Dehradun, India}

\email{rahulkarki877@gmail.com} 
\thanks{}

\author{B. Subhash}

\address{Indian Institute of Science Education and Research, Tirupati, Andhra Pradesh, India.} 

\email{subhash@iisertirupati.ac.in}
\thanks{}


\begin{document}

\begin{abstract}
    The \textit{perfect matching complex} of a simple graph $G$ is an abstract simplicial complex whose facets (maximal faces) are precisely the perfect matchings of $G$. 
    In this article, we study the topology of the perfect matching complex of the $3\times 2n$ grid graph  and show that it is homotopy equivalent to a wedge of $f_{n-1}$ copies of $(n-1)$-dimensional spheres, where $f_{k}$ denotes the $k$-th Fibonacci number. This proves a particular case of the conjecture given by Chandrakar and Singh for perfect matching complex of grid graphs.\\
    % Interestingly, the number of spheres appearing in this homotopy equivalence is given by a Fibonacci number. More precisely, we show that the perfect matching complex is homotopy equivalent to a wedge of $f_{n-1}$ spheres of dimension $n-1$, where $f_{n-1}$ denotes the $(n-1)$-th Fibonacci number.
\end{abstract}

\maketitle    

\section{Introduction}

Matchings in a graph are fundamental objects in graph theory and have been studied extensively from both combinatorial and algebraic points of view (\cite{plummer1986matching}). The collections of matchings of a graph naturally give rise to a simplicial complex, called the \textit{matching complex}. Over the years, the topology of matching complexes has been studied extensively, leading to several interesting results concerning their homotopy types and connectivity (\cite{christos, bjorner1994chessboard, jonsson2008simplicial, matsushita2022matching, Matsushita2018MatchingCO, wachs2003topology}). It is also worth mentioning that the occurrence of torsion in the homology groups of matching complexes highlights intriguing and subtle aspects of their topology (\cite{bouc1992homologie, wachs2003topology, shareshian2007torsion}).

Another interesting graph theoretic concept closely related to matchings is that of perfect matchings. Recently, Bayer et al. introduced the notion of the \textit{perfect matching complex} of a graph, which is a subcomplex of its matching complex \cite{bayer2022perfect}. Since its introduction, the topology of perfect matching complexes has attracted considerable attention.

Bayer et al. initiated the study of perfect matching complexes of honeycomb graphs (or hexagonal grid graphs), with the broader aim of understanding the topology of matching complexes associated with these graphs. In particular, they studied the perfect matching complexes of hexagonal grid graphs. Their approach combines discrete Morse theory with combinatorial methods involving plane partitions and other topological techniques.

% \sgnote{We should mention that Chandrakar et al. conjectured that the  perfect matching complex of grid graphs is homotopy equivalent to a wedgeof spheres. We settle this conjecture for the case $3 \times n$. Maybe we should mention it in the abstract also?}
Subsequently, Chandrakar et al. studied perfect matching complexes of polygonal line tilings, with particular emphasis on the \(2\times n\) grid graph, also known as the \textit{ladder graph} \cite{ChandrakarSingh2025}. 
They introduced the notion of \textit{bad matchings} (non-extendable matchings), namely, matchings that cannot be extended to perfect matchings, and used this characterization together with discrete Morse theory and the fold lemma for independence complexes (\cite{engstrom_clawfree, forman2002user}) to determine the homotopy types of these complexes. 
More recently, Agarwal et al. investigated perfect matching complexes of M{\"o}bius ladder graphs \cite{agarwal2026_pmc_laddergraph}.
Moreover, the matching complexes of $3\times n$ grid graphs have recently been studied by Chauhan et al. \cite{chauhan2026_matching3n}, who determined their topological connectivity.

The approach based on bad matchings is particularly effective for the $2\times n$ grid graph. A bad matching (see \cite[Definition~2.3]{ChandrakarSingh2025}) is a matching that cannot be extended to a perfect matching but whose proper subsets can all be extended to perfect matchings. Thus, bad matchings correspond precisely to the minimal nonfaces of the perfect matching complex. Chandrakar et al. showed that every bad matching in the $2\times n$ grid graph consists of exactly two edges (see \cite[Lemma~3.1]{ChandrakarSingh2025}). Recall that a simplicial complex is called a \textit{flag complex} if every minimal nonface has exactly two elements. Moreover, any two edges that share a vertex cannot belong to the same matching and hence form a minimal nonface of the perfect matching complex. Therefore, every minimal nonface of the perfect matching complex of the $2 \times n$ grid graph has size two, and hence the complex is a flag complex. Consequently, it can be realized as the independence complex of a suitable graph, making results from the theory of independence complexes available for studying its topology.

The next natural case to consider in the family of $m\times n$ grid graphs is the $3\times n$ grid graph. Since each edge covers exactly two vertices, a graph admitting a perfect matching must have an even number of vertices. Thus, for the $3\times n$ grid graph to admit a perfect matching, $n$ must be even. Accordingly, in this article, we study the perfect matching complex of the $3\times 2n$ grid graph. However, extending the approach based on bad matchings from the $2\times n$ grid graph to this family is not straightforward. While every bad matching in the former consists of exactly two edges, bad matchings consisting of three edges can occur in the latter. Consequently, the perfect matching complex of the $3\times 2n$ grid graph need not be a flag complex and cannot, in general, be realized as the independence complex of a graph. Thus, the methods based on this correspondence do not directly extend to this case, making the computation of its homotopy type more challenging.

In \cite{ChandrakarSingh2025}, it is conjectured that the perfect matching complex of grid graphs is homotopy equivalent to a wedge of spheres. In this article, we settle this conjecture for the case $3 \times 2n$. 
% Interestingly, the number of spheres appearing in its homotopy type is governed by the Fibonacci sequence. 
More precisely, we prove the following result.

\begin{theorem}[\Cref{thm: the main theorem}]
    The perfect matching complex of the $3\times 2n$ grid graph is contractible for $n=1$, while for $n\geq 2$, it is homotopy equivalent to a wedge of $f_{n-1}$ spheres of dimension $n-1$, where $f_k$ denotes the $k$th Fibonacci number, with $f_1=f_2=1$.
\end{theorem}

The article is organized as follows. Section 2 presents the necessary preliminaries, including the relevant definitions and results used throughout the article. In Section 3, we begin by formally defining the $3 \times 2n$ grid graph and specifying the labelling of its vertices and edges. We then establish a few lemmas, %the main lemma, 
which are used to prove our main theorem.

\section{Preliminaries}

Let $G$ be an undirected simple graph with vertex set $V(G)$ and edge set $E(G)$. Then a subset $M$ of $E(G)$ is said to be a {\it matching} if endpoints of no two edges in $M$ are same. Furthermore, a matching $M$ is said to be a {\it perfect matching} if the set of end points of all the edges of $M$ equals $V(G)$. 
% \sgnote{put a small example}

\begin{example}
    Let $G$ be the 6-cycle graph with $V(G)=\{x_1,\dots, x_6\}$ and $E(G)=\{\{x_i,x_{i+1}\} \mid i=1,\dots, 6\}$, where $x_7\coloneqq x_1$. The set of edges $\{ \{x_1,x_2\},\{x_3,x_4\}, \{x_5,x_6\} \}$ is a perfect matching since all vertices appear as endpoints. On the other hand, the matching $\{  \{x_1,x_2\}, \{x_4,x_5\}\}$ is not perfect, despite it being a maximal matching, in the sense that no strictly larger matching contains it.
\end{example}

% \subsection{Simplicial Complexes}

We begin by recalling some standard definitions concerning simplicial complexes.

\medskip

Let $V$ be a finite set. An {\it (abstract) simplicial complex} $K$ on the set $V$ is a family of subsets of $V$ with the property that whenever $F \in K$, every subset of $F$ is also contained in $K$. The members of $K$ are called {\it faces}, and the inclusion-maximal faces are called {\it facets}. 
%We denote the set of facets of $\Delta$ by $\F(\Delta)$.

If $F = \{x_1, x_2, \dots, x_n\} \subseteq V$ is a facet, we denote by $\langle x_1, x_2, \dots, x_n \rangle$ the subcomplex of $K$ given by the power set of $F$. 

The {\it dimension} of a face $F$ is given by $\dim F = |F|-1$, where $|F|$ denotes the number of elements in $F$, and in such a case $F$ is said to be $(|F|-1)$-dimensional face of $K$. A 0-dimensional face of $K$ is called a {\it vertex} of $K$. The {\it dimension} of the simplicial complex $K$ is defined as the largest dimension among its faces, i.e., 
\[\dim(K) = \max_{F \in K}\{\dim(F)\}.\]

Let $K$ and $L$ be two simplicial complexes, then {\it join of simplicial complexes} $K$ and $L$ is defined as
\[K \ast L = \{F_K \cup F_L| F_K \in K, F_L \in L\}.\]

Let $F \in K$ be a face.
The {\it link}, $\lk_K(F)$ and {\it deletion}, $\del_K(F)$ of $F$ in $K$ are subcomplexes of $K$ defined as
\[\lk_K(F) = \{F' \in K : F \cap F' = \emptyset, F \cup F' \in K\},\;
\del_K(F) = \{F' \in K : F' \cap F = \emptyset\}.\]

Let $a,b$ be two vertices that are not in $K$. The {\it suspension} $\Sigma K$ of $K$ is defined as
\[
\Sigma K = K\cup \{\sigma\cup \{a\}, \sigma\cup \{b\}\mid \sigma\in K\}.
\]

The {\it wedge} $K\vee L$ is the simplicial complex obtained by identifying one vertex each from $K$ and $L$.
% \sgnote{we need to be more careful for the choice of vertex, for $K$ or $L$ can have more than one connected components}

	\begin{lemma}[{\cite[Proposition~3.1]{ADAMASZEK20121031}}] \label{lemma: link-deletion}
		Let $K$ be a simplicial complex and $v$ be a vertex of $K$. If the inclusion map $\lk_K(v) \to \del_K(v)$ is null-homotopic, then 
		$$ K \simeq \del_K(v) \vee \Sigma (\lk_K(v)),$$
		where $\vee$ denotes the wedge of spaces, and $\Sigma$ denotes the suspension of space.
	\end{lemma}

    % \sgnote{define wedge and suspension?}

The {\it star of a face $F \in K$}, $\st_K(F)$ is a subcomplex of $K$ defined as 
\[\st_K(F) = \{F' \in K : F \cup F' \in K\}.\]
It is easy to verify that $\st_K(F)$ is a contractible subcomplex of $K$ for any face $F$ of $K$.

\begin{lemma}[{\cite[Lemma~10.4]{BjornerTopoMethods}}]\label{lemma: Bjorner-star-intersection}
    Let $\Delta = \Delta_0 \cup \Delta_1 \cup \dots \cup \Delta_n$ be a simplicial complex with contractible subcomplexes $\Delta_i$ such that $\Delta_i \cap \Delta_j \subseteq \Delta_0$ for $0< i<j \le n$, then 
    \[\Delta \simeq \bigvee_{1 \le i \le n} \Sigma(\Delta_o \cap \Delta_i).\]
\end{lemma}
In this article, we will use the particular case of the above theorem for $n=1$.

\section{Main results}

Let $m,n$ be natural numbers. The $m \times n$ grid graph $\mathcal{G}_{m \times n}$ is a graph with the vertex set $V(\mathcal{G}_{m \times n}) = \left\{ p_{i,j}  :  1 \le i \le m, 1 \le j \le n\right\}$ and edges $E(\mathcal{G}_{m \times n}) = \{\{p_{i,j}, p_{k,l}\}  : |i-k| + |j-l|=1 \}$.

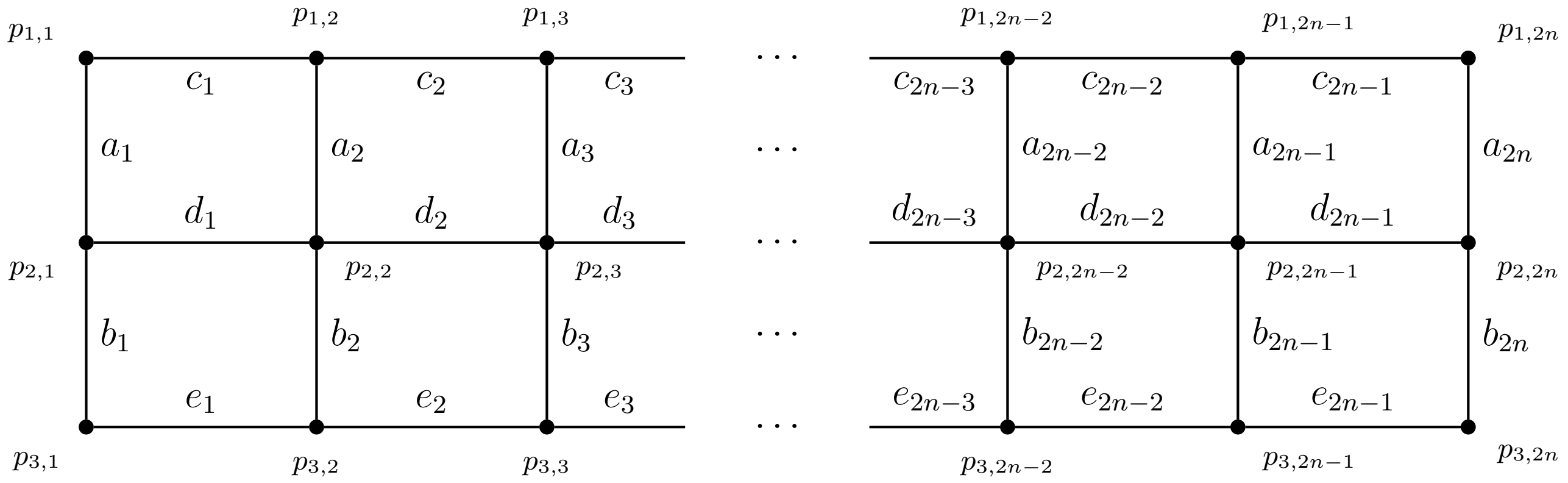
\begin{figure}[htbp]
    \centering
    \begin{tikzpicture}[
        vertex/.style={circle, draw=black, fill=black, inner sep=1.5pt},
        edge/.style={thick},
        every label/.style={font=\scriptsize},
        label distance=1mm
    ]

    % 1. VERTICES
    % Top Row (i = 1)
    \node[vertex,label={175:$p_{1,1}$}]  (T1)  at (0, 4) {};
    \node[vertex,label={90:$p_{1,2}$}]   (T2)  at (2.5, 4) {};
    \node[vertex,label={90:$p_{1,3}$}]   (T3)  at (5, 4) {};
    \node[vertex,label={90:$p_{1,2n-2}$}] (Tn1) at (10, 4) {};
    \node[vertex,label={45:$p_{1,2n-1}$}]  (Tn2)  at (12.5, 4) {};
    \node[vertex,label={5:$p_{1,2n}$}]  (Tn3)  at (15, 4) {};

    % Middle Row (i = 2)
    \node[vertex,label={200:$p_{2,1}$}]  (M1)  at (0, 2) {};
    \node[vertex,label={340:$p_{2,2}$}]  (M2)  at (2.5, 2) {};
    \node[vertex,label={340:$p_{2,3}$}]  (M3)  at (5, 2) {};
    \node[vertex,label={340:$p_{2,2n-2}$}] (Mn1) at (10, 2) {};
    \node[vertex,label={340:$p_{2,2n-1}$}] (Mn2)  at (12.5, 2) {};
  \node[vertex,label={340:$p_{2,2n}$}] (Mn3)  at (15, 2) {};

    % Bottom Row (i = 3)
    \node[vertex,label={225:$p_{3,1}$}]  (B1)  at (0, 0) {};
    \node[vertex,label={270:$p_{3,2}$}]  (B2)  at (2.5, 0) {};
    \node[vertex,label={270:$p_{3,3}$}]  (B3)  at (5, 0) {};
    \node[vertex,label={270:$p_{3,2n-2}$}] (Bn1) at (10, 0) {};
    \node[vertex,label={315:$p_{3,2n-1}$}] (Bn2)  at (12.5, 0) {};
   \node[vertex,label={345:$p_{3,2n}$}] (Bn3)  at (15, 0) {};

    % 2. HORIZONTAL EDGES
    \draw[edge] (T1) -- node[below] {$c_1$} (T2);
    \draw[edge] (T2) -- node[below] {$c_2$} (T3);
    \draw[edge] (T3) -- node[below] {$c_3$} (6.5, 4);
    \path (6.5, 4) -- (8.5, 4) node[midway] {$\cdots$};
    \draw[edge] (8.5, 4) -- node[below] {$c_{2n-3}$} (Tn1);
    \draw[edge] (Tn1) -- node[below] {$c_{2n-2}$} (Tn2);
    \draw[edge] (Tn2) -- node[below] {$c_{2n-1}$} (Tn3);

    \draw[edge] (M1) -- node[above] {$d_1$} (M2);
    \draw[edge] (M2) -- node[above] {$d_2$} (M3);
    \draw[edge] (M3) -- node[above] {$d_3$} (6.5, 2);
    \path (6.5, 2) -- (8.5, 2) node[midway] {$\cdots$};
    \draw[edge] (8.5, 2) -- node[above] {$d_{2n-3}$} (Mn1);
    \draw[edge] (Mn1) -- node[above] {$d_{2n-2}$} (Mn2);
    \draw[edge] (Mn2) -- node[above] {$d_{2n-1}$} (Mn3);

    \draw[edge] (B1) -- node[above] {$e_1$} (B2);
    \draw[edge] (B2) -- node[above] {$e_2$} (B3);
    \draw[edge] (B3) -- node[above] {$e_3$} (6.5, 0);
    \path (6.5, 0) -- (8.5, 0) node[midway] {$\cdots$};
    \draw[edge] (8.5, 0) -- node[above] {$e_{2n-3}$} (Bn1);
    \draw[edge] (Bn1) -- node[above] {$e_{2n-2}$} (Bn2);
    \draw[edge] (Bn2) -- node[above] {$e_{2n-1}$} (Bn3);

    % 3. VERTICAL EDGES
    \draw[edge] (T1) -- node[right] {$a_1$} (M1);
    \draw[edge] (T2) -- node[right] {$a_2$} (M2);
    \draw[edge] (T3) -- node[right] {$a_3$} (M3);
    \draw[edge] (Tn1) -- node[right] {$a_{2n-2}$} (Mn1);
    \draw[edge] (Tn2) -- node[right] {$a_{2n-1}$} (Mn2);
    \draw[edge] (Tn3) -- node[right] {$a_{2n}$} (Mn3);

    \draw[edge] (M1) -- node[right] {$b_1$} (B1);
    \draw[edge] (M2) -- node[right] {$b_2$} (B2);
    \draw[edge] (M3) -- node[right] {$b_3$} (B3);
    \draw[edge] (Mn1) -- node[right] {$b_{2n-2}$} (Bn1);
    \draw[edge] (Mn2) -- node[right] {$b_{2n-1}$} (Bn2);
    \draw[edge] (Mn3) -- node[right] {$b_{2n}$} (Bn3);

    % 4. GRID CONTINUATION DOTS
    \path (6.5, 3) -- (8.5, 3) node[midway] {$\cdots$};
    \path (6.5, 1) -- (8.5, 1) node[midway] {$\cdots$};

    \end{tikzpicture}
    \caption{The grid graph $\mathcal{G}_{3 \times 2n}$}
    \label{grid graph}
\end{figure}

For a graph $G$, let $\mathcal{M}_p(G)$ denote the {\it perfect matching complex} of $G$, i.e., the simplicial complex with facets being the perfect matchings of $G$.\\
% Let $\mathcal{G}_{m,n}$ be the $m\times n$ grid graph.
We denote the perfect matching complex $\pmcggg$ by $\Delta_n.$ 
For $1 \le j \le 2n$, let $a_j, b_j$ denote the edges $\{p_{1,j}, p_{2,j}\}$, $\{p_{2,j}, p_{3,j}\}$, respectively and for $1 \le j \le 2n-1$, let $c_j = \{p_{1,j}, p_{1,j+1}\}$, $d_j = \{p_{2,j}, p_{2,j+1}\}$ and $e_j =\{p_{3,j}, p_{3,j+1}\}$, see \Cref{grid graph}.

\pagebreak

%%% ∆n is union of star(c(2n-1)) and sta(a(2n)) %%%
\begin{lemma}\label{lemma: Delta_n is union of stars}
    For $n \geq 1$, we have $$\Delta_n = \st_{\Delta_n}(a_{2n}) \cup \st_{\Delta_n}(c_{2n-1}).$$
\end{lemma}

\begin{proof}
    Note that any perfect matching of $\mathcal{G}_{3\times 2n}$ should contain either $c_{2n-1}$ or $a_{2n}$ to cover the vertex $p_{1,2n}.$ Thus, every facet of $\Delta_n = \pmcggg$ will contain either $c_{2n-1}$ or $a_{2n}.$ The result then follows from the definition of the star of a vertex.
\end{proof}

%\sgnote{}

For the ease of notation, we will write $\st(v)$ to denote the $\st_{\Delta_n}(v)$. For $n\ge 1$, let us denote the intersection of $\st(a_{2n})$ and $\st(c_{2n-1})$ by $X_n$, {\it i.e.,}
\[X_n= \st(a_{2n}) \cap \st(c_{2n-1}).\]

 Let $\sigma$ be a facet of $\st( a_{2n})$, then the edges $b_{2n}$, $c_{2n-1}$, and $d_{2n-1}$ being adjacent to $a_{2n}$, we get $\sigma \cap \{b_{2n}, c_{2n-1}, d_{2n-1}\} = \emptyset$. 
    Moreover, to cover the vertex $p_{3,2n}$, the edge  $e_{2n-1}$ should belong to $\sigma$, implying that $e_{2n-2}$ and $b_{2n-1}$ cannot belong to $\sigma$. 
    This implies that  $ \sigma \cap \{b_{2n}, b_{2n-1}, c_{2n-1}, d_{2n-1}, e_{2n-2}\} = \emptyset.$

   \vspace{2mm}
    Similarly, if $\sigma$ is a facet of $\st_{\Delta_n} \left( c_{2n-1} \right)$, then the edges $a_{2n}$, $a_{2n-1}$, and $c_{2n-2}$ will not be in $\sigma$ as they are adjacent to $c_{2n-1}$, that is,  $\sigma \cap \{a_{2n}, a_{2n-1}, c_{2n-2}\} = \emptyset.$

    \noindent Thus, if $\sigma \in \st \left( a_{2n} \right) \cap \st \left( c_{2n-1} \right)$, then 
    \begin{equation}\label{eqn: edges that sigma in Yn not contain}
        \sigma \cap \{a_{2n}, a_{2n-1}, b_{2n}, b_{2n-1}, c_{2n-1}, c_{2n-2}, d_{2n-1}, e_{2n-2}\} = \emptyset.
    \end{equation}

  If $\tau$ is a facet of $X_n,$ then since $\tau \in \st \left( c_{2n-1} \right),$ $\tau \cup \left\{ c_{2n-1} \right\} \in \Delta_n$. 
   Let $\tau_{c_{2n-1}}$ be a facet of $\Delta_n = \pmcggg$ containing $\tau \cup \left\{ c_{2n-1} \right\},$ then $\tau_{c_{2n-1}}$ is a perfect matching containing $c_{2n-1},$ extending the facet $\tau \in X_n.$ Similarly, let $\tau_{a_{2n}}$ be a perfect matching of $\Delta_n = \pmcggg$ containing $a_{2n},$ extending the facet $\tau \in X_n.$ 

   \begin{lemma} \label{lemma: edges d2n-2 e2n-1}
     If $\tau$ is a facet of $X_n$ then $d_{2n-2}$ and $e_{2n-1}$ cannot appear simultaneously in $\tau,$ that is $\{d_{2n-2}, e_{2n-1}\} \nsubseteq \tau.$
   \end{lemma}
   \begin{proof}
     On the contrary, if $\left\{d_{2n-2}, e_{2n-1} \right\} \subseteq \tau$, then $\tau_{c_{2n-1}}$ is a perfect matching containing $\{c_{2n-1},d_{2n-2}, e_{2n-1}\}.$ 
    % Let $\tilde{\tau} \in \Delta_n$ be a facet containing $\tau \cup \{c_{2n-1}\}$, 
     Hence $\tau_{c_{2n-1}} \cap \{d_{2n-1}, a_{2n}, b_{2n}\} = \emptyset,$ which implies that the vertex $p_{2,2n}$ is uncovered, contradicting the fact that $\tau_{c_{2n-1}}$ is a perfect matching.
   \end{proof}
% \sgnote{ The defnition $A_k$ below requires a little more explanation}

   \noindent For $1 \le k \le n-1$, let $\mathcal{A}_k$ be the subcomplex of $\Delta_n$ given by  \[\mathcal{A}_k = \begin{cases} \langle d_{2n-2}, d_{2n-4}, \dots, d_4, d_2 \rangle  \qquad \qquad \qquad \mbox{ if } k=1,\\\langle d_{2n-2}, d_{2n-4}, \dots, d_{2k+2}, d_{2k} \rangle * \Delta_{k-1} \quad \, \mbox{ if } 2 \leq k \leq n-1    .\end{cases}\]

   \begin{lemma} \label{Lemma : facets containing d_{2n-2}}
     If $\tau $ is a facet of $X_n$ and if $d_{2n-2} \in \tau,$ then $\tau \in \bigcup\limits_{i=1}^{n-1} \mathcal{A}_{i}$
   \end{lemma}
   \begin{proof}
     Since $d_{2n-2} \in \tau,$ we have
     \begin{equation}\label{eqn: edges not in tau due to d_2n-2}
        \tau \cap \{a_{2n-1}, a_{2n-2}, b_{2n-1}, b_{2n-2}, d_{2n-3}, d_{2n-1}\} = \emptyset.
     \end{equation}
     By \Cref{lemma: edges d2n-2 e2n-1}, we also have $e_{2n-1} \notin \tau.$
 
   \noindent For $2 \leq i \leq n$, we define 
    $$T_i = \begin{cases}
        \left\{ a_{2i-2},\ a_{2i-3},\ b_{2i-2},\ b_{2i-3},\ c_{2i-3},\  d_{2i-3},\ e_{2i-3}\right\} & \text{ for } i=2,\\
        \left\{ a_{2i-2},\ a_{2i-3},\ b_{2i-2},\ b_{2i-3},\ c_{2i-4}, \ c_{2i-3},\  d_{2i-3},\ e_{2i-4}, \ e_{2i-3}\right\} & \text{ otherwise}.
    \end{cases}$$

    \noindent We claim that $\tau \cap T_n = \emptyset.$   
    Since $\tau_{c_{2n-1}}$ is a perfect matching containing $d_{2n-2},$ in order to cover the vertex $p_{1,2n-2}$ it should contain $c_{2n-3}$,  and thus $a_{2n-3}, c_{2n-4} \notin \tau_{c_{2n-1}},$ therefore not in $\tau.$ 
% hence not in $\tau$.
    Similarly, since $\tau_{a_{2n}}$ is a perfect matching containing $d_{2n-2},$ to cover the vertex $p_{1,2n-1}$ it should contain $c_{2n-2},$   hence $c_{2n-3} \notin \tau_{a_{2n}},$ and thus not in $\tau$ either.
    
    \noindent In order to cover the vertex $p_{3,2n-2}$, the perfect matching $\tau_{c_{2n-1}}$ should contain $e_{2n-2}$, hence $e_{2n-3} \notin \tau_{c_{2n-1}}$ and the perfect matching $\tau_{a_{2n}}$ should contain $e_{2n-3}$, hence $e_{2n-4}, b_{2n-3} \notin \tau_{a_{2n}}$. 
    This implies $e_{2n-3}, e_{2n-4}, b_{2n-3} \notin \tau,$ hence we have $\tau \cap T_n = \emptyset.$

     We now examine the remaining \(d\)-edges. 
    Since \(d_{2n-2}\in\tau\) and \(\tau\cap T_n=\emptyset\), the next \(d\)-edge to be considered is \(d_{2n-4}\). If \(d_{2n-4} \notin \tau\), then $\tau \in \langle d_{2n-2}\rangle \ast \Delta_{n-1} = \mathcal{A}_{n-1}.$
    
    \noindent If $d_{2n-4} \in \tau$, then the above argument can be repeated with the indices shifted by two. 
    Thus, $\tau \cap T_{n-1} = \emptyset$. Since $\tau \cap T_n = \emptyset$, we get $\tau \cap (T_n \cup T_{n-1}) = \emptyset$. 
    Moreover, as before, if $d_{2n-6} \notin \tau$, then $\tau \in \langle d_{2n-2}, d_{2n-4} \rangle \ast \Delta_{n-2} = \mathcal{A}_{n-2}$. 

   % Proceeding similarly for the remaining $d$-edges, let $1 \leq k \leq n-1$ be the largest $k$ such that for every $i \leq k$, $d_{2n-2i} \in \tau$ but $d_{2n-2k-2} \notin \tau$.

    \noindent Proceeding similarly, if 
    $d_{2n-2}, d_{2n-4} \dots, d_{2n-2k+2}, d_{2n-2k} \in \tau$ where $1 \leq k \leq n-1,$ then 
    $$\tau \cap \left( \bigcup\limits_{i=0}^{k-1}T_{n-i} \right) = \emptyset.$$ 
    If $d_{2n-2(k-1)} \notin \tau$, then $\tau \in \mathcal{A}_{n-k}.$ This shows that $\tau \in  \bigcup\limits_{i=1}^{n-1} \mathcal{A}_{i}.$
\end{proof}

    \noindent With the two lemmas above, we can prove the following result.

\begin{lemma}\label{lemma: intersection of stars}
    For $n\geq 2$, 
    we have,
    \[
        X_n = \left( \bigcup\limits_{i=1}^{n-1} \mathcal{A}_{i}\right) \cup \left(\langle e_{2n-1}  \rangle * \Delta_{n-1}\right),
    \]  where $\mathcal{A}_k = \begin{cases} \langle d_{2n-2}, d_{2n-4}, \dots, d_4, d_2 \rangle  \quad \qquad \qquad \mbox{ if } k=1,\\\langle d_{2n-2}, d_{2n-4}, \dots, d_{2k+2}, d_{2k} \rangle * \Delta_{k-1} \,\mbox{ if } 2 \leq k \leq n-1    .\end{cases}$
\end{lemma}

\begin{proof}
    We first show that, $X_n \subseteq \left( \bigcup\limits_{i=1}^{n-1} \mathcal{A}_{i}\right) \cup \left(\langle e_{2n-1} \rangle * \Delta_{n-1} \right).$ 
    Let $\tau$ be a facet of $X_n,$ by \Cref{lemma: edges d2n-2 e2n-1}, since $e_{2n-1}$ and $d_{2n-2}$ cannot simultaneously belong to $\tau.$
    We have
    \vspace{2mm}
    
  \noindent  \textbf{Case 1:} Consider the case when $e_{2n-1} \in \tau$.
  
    \noindent Observe that for $\tau_{c_{2n-1}},$ we have $c_{2n-1},e_{2n-1} \in \tau_{c_{2n-1}}$. 
    In order to cover $p_{2,2n},$ it follows that $d_{2n-1} \in \tau_{c_{2n-1}},$ hence  $\tau_{c_{2n-1}} \cap \left\{ d_{2n-2}\right\} = \emptyset$. 

    \noindent This along with {\Cref{eqn: edges that sigma in Yn not contain}},  implies that $$\tau\cap\left\{a_{2n},b_{2n}, c_{2n-1}, a_{2n-1}, b_{2n-1}, d_{2n-1}, d_{2n-2}, c_{2n-2}, e_{2n-2} \right\} = \emptyset.$$

    \noindent Hence we have $\tau \in \langle  e_{2n-1} \rangle * \Delta_{n-1}.$

  %  \noindent This completes case 1.
   \vspace{2mm}
 
    \noindent \textbf{Case 2:} Suppose $e_{2n-1} \notin \tau.$ 
    Since every facet in $\st(a_{2n})$ contains $e_{2n-1},$ we have $e_{2n-1} \notin \tau_{c_{2n-1}}.$ 
    If $d_{2n-2} \notin \tau,$ {\Cref{eqn: edges that sigma in Yn not contain}} implies $\tau \in \Delta_{n-1},$ in which case $\tau \cup \{e_{2n-1}\} \in X_n,$ implies that $\tau$ is not a facet of $X_n.$ Hence $d_{2n-2} \in \tau.$ 
    Then by  \Cref{Lemma : facets containing d_{2n-2}}, we have $\tau \in\bigcup\limits_{i=1}^{n-1} \mathcal{A}_{i}. $\\
   This shows that $X_n \subseteq \left( \bigcup\limits_{i=1}^{n-1} \mathcal{A}_{i}\right) \cup \left(\langle e_{2n-1} \rangle * \Delta_{n-1} \right).$

    \vspace{3mm}
    %%% Reverse Inclusion %%%

    \noindent We now prove the reverse inclusion.
    % that is, $$X_n \supseteq \left( \bigcup\limits_{i=1}^{n-1} \mathcal{A}_{i}\right) \cup \left(\langle \left( e_{2n-1} \right) \rangle * \Delta_{n-1} \right).$$
    Let $\tau \in \left( \bigcup\limits_{i=1}^{n-1} \mathcal{A}_{i}\right) \cup \left(\langle e_{2n-1}\rangle * \Delta_{n-1} \right)$ be a facet. We need to show that $\tau \in X_n = \st (c_{2n-1}) \cap \st (a_{2n} )$.

   If $\tau \in \langle e_{2n-1} \rangle * \Delta_{n-1}$, then $\tau \cup \left\{c_{2n-1}, d_{2n-1}\right\}$ is a perfect matching containing $c_{2n-1}$ and hence belongs to $\st({c_{2n-1}})$. Similarly, $\tau\cup \left\{a_{2n-1}, a_{2n}\right\}$ is a perfect matching that contains $a_{2n}$ and hence belongs to $\st(a_{2n})$. This implies $\tau \in X_n$. 

    \noindent If $\tau \in  \bigcup\limits_{i=1}^{n-1} \mathcal{A}_{i} ,$ then $\tau_1$ defined as $$\tau_1 = \begin{cases} \tau \cup \left\{ b_1, b_{2n}, c_1, c_3, \dots, c_{2n-1}, e_2, e_4, \dots, e_{2n-2} \right\} \mbox { if } \tau  \in \mathcal{A}_1 = \langle d_{2n-2}, d_{2n-4}, \dots, d_{2} \rangle \\ \tau \cup \{ b_{2i-1}, b_{2n}, c_{2i-1}, c_{2i+1}, \dots, c_{2n-1}, e_{2i}, e_{2i+2}, \dots, e_{2n}\} \mbox{ if } \tau \in \mathcal{A}_i, 2\leq i \leq n-1, \end{cases}$$ 
    is a perfect matching containing $c_{2n-1}$ and hence belongs to $\st({c_{2n-1}})$.
    Similarly $\tau_2$ defined as $$ \tau_2 = \begin{cases}\tau \cup \left\{ a_1, a_{2n}, c_2, c_4, \dots, c_{2n-2}, e_1, e_3, \dots, e_{2n-1} \right\} \mbox { if } \tau  \in \mathcal{A}_1 = \langle d_{2n-2}, d_{2n-4}, \dots, d_{2} \rangle \\\tau \cup \{ a_{2i-1}, a_{2n}, c_{2i}, c_{2i+2}, \dots, c_{2n-2}, e_{2i-1}, e_{2i+1}, \dots, e_{2n-1}\} \mbox{ if } \tau \in \mathcal{A}_i, 2\leq i \leq n-1,  \end{cases}$$ is a perfect matching containing ${a_{2n}}$, therefore $\tau_2 \in \st({a_{2n}})$ and hence $\tau \in X_n$.
    \noindent From the above discussion, we have $\left( \bigcup\limits_{i=1}^{n-1} \mathcal{A}_{i}\right) \cup \left(\langle e_{2n-1}\rangle * \Delta_{n-1} \right) \subseteq X_n$.
    This completes the proof.
\end{proof} 

%%% Main theorem %%%
\begin{theorem}\label{thm: the main theorem}
    Let $\Delta_n$ be the perfect matching complex of $3 \times 2n$ grid graph, then $\Delta_1$ is contractible, and for $n \geq 2$, 
    \[\Delta_n \cong \displaystyle{\bigvee_{f_{n-1}} \mathbb{S}^{n-1}}, 
    \]
    where $f_k$ denotes the $k$-th Fibonacci number with $f_1 = f_2 = 1$.\\ 
    % \sgnote{should we include $\Delta_1$ also in the statement here itself}
\end{theorem}

\begin{proof}
    Note that since $\Delta_1 = \langle a_1, a_2, e_1\rangle \cup \langle b_1, b_2, c_1 \rangle \cup \langle c_1, d_1, e_1 \rangle $, it is easy to see that $\Delta_1$ is contractible. 
    For $n \geq 2$, from \Cref{lemma: Delta_n is union of stars}, we get $\Delta_n = \st(a_{2n}) \cup \st(c_{2n-1}).$
    Since both $\st(a_{2n})$ and $ \st(c_{2n-1})$ are contractible, we apply \Cref{lemma: Bjorner-star-intersection} to get \begin{equation}\label{eqn: homotopy type of Delta-n}
        \Delta_n = \Sigma \left( \st(a_{2n}) \cap \st(c_{2n-1})\right).
    \end{equation}

   \noindent Therefore, it suffices to determine the homotopy type $X_n= \st \left( a_{2n} \right) \cap \st \left( c_{2n-1} \right)$. 
   From \Cref{lemma: intersection of stars}, we have \begin{equation}\label{eqn: short expression of Xn}
        X_n = \left( \bigcup\limits_{i=1}^{n-1} \mathcal{A}_{i}\right) \cup \left(\langle e_{2n-1} \rangle * \Delta_{n-1} \right),
    \end{equation}
   where, $$\mathcal{A}_i = \begin{cases}
        \langle d_{2}, d_{4}, \dots, d_{2n-2} \rangle &\mbox{ if } i=1,\\
        \langle d_{2i}, d_{2i+2}, \dots, d_{2n-2} \rangle * \Delta_{i-1} & \mbox{ if } 2 \leq i \leq n-1.
    \end{cases}$$

  %\noindent We determine the homotopy type of $\Delta_2$. 
  \noindent Using \Cref{eqn: short expression of Xn}, we obtain $X_2 = \left \langle d_2 \right \rangle \cup \left( \left \langle e_3  \right \rangle * \Delta_1 \right).$
  Since both components are contractible and disjoint, it follows that $ X_2 \simeq \mathbb{S}^0.$
  Hence, \begin{equation*}\label{eqn: homotopy of Delta_2}
        \Delta_2 = \Sigma X_2 \simeq \mathbb{S}^1.
    \end{equation*}

    %%% deletion of d(2n-2) in Xn %%%

    \noindent Let $n \geq 3$ and $\tau \in X_n$ be a facet. It is easy to see that if $d_{2n-2} \notin \tau$, then $\tau \in \del_{X_n}(d_{2n-2})$, else $\tau \setminus \{d_{2n-2}\} \in \del_{X_n}(d_{2n-2})$ is a facet. This along with \Cref{eqn: short expression of Xn} implies 

     $$\del_{X_n} \left( d_{2n-2} \right) = \left( \bigcup\limits_{i=1}^{n-1} \mathcal{A}_{i}' \right) \cup \left(\langle  e_{2n-1} \rangle * \Delta_{n-1} \right),$$

    \noindent where $\mathcal A_i'$ is obtained by removing $d_{2n-2}$ from every facet of $\mathcal A_i$. That is, $$\mathcal{A}_i' = \begin{cases}
        \langle  d_{2}, d_{4}, \dots, d_{2n-4} \rangle & \mbox{ if }i=1,\\
        \langle  d_{2i}, d_{2i+2}, \dots, d_{2n-4} \rangle * \Delta_{i-1} &\mbox { if } 2 \leq i \leq n-2,\\
        \Delta_{n-2} & \mbox{ if }i = n-1.
    \end{cases}$$

    \noindent Since $d_i \in \Delta_{n-1}$ for $i \leq 2n-3$, it follows that  $\mathcal{A}_k'$ is a subcomplex of $\Delta_{n-1}$ for $1 \leq k \leq n-1$. 
    Therefore, $\del_{X_n}(d_{2n-2}) = \langle e_{2n-1} \rangle * \Delta_{n-1}$, implying that $\del_{X_n}(d_{2n-2})$ is contractible.
    Hence, the inclusion $\lk_{X_n} \left( d_{2n-2} \right) \hookrightarrow\del_{X_n} \left( d_{2n-2} \right)$ is null-homotopic. 
    Therefore by \Cref{lemma: link-deletion}, $X_n \simeq \del_{X_n} \left( d_{2n-2} \right) \vee \Sigma \left( \lk_{X_n} \left( d_{2n-2} \right) \right) \simeq \Sigma ( \lk_{X_n} ( d_{2n-2}))$, and hence $\Delta_{n} \simeq \Sigma^2 ( \lk_{X_n} ( d_{2n-2} )).$ Therefore it suffices to determine the homotopy type of
    $\lk_{X_n}(d_{2n-2})$.
        
    \noindent For $n \geq 3$, let $Z_n=\lk_{X_n}(d_{2n-2})$.
    From \Cref{eqn: short expression of Xn}, 
    we see that 
    \begin{equation}\label{eqn: descriptive form of Zn}
        Z_n = \lk_{X_n}( d_{2n-2}) =  \bigcup\limits_{i=1}^{n-1} \mathcal{A}_{i}' = \Delta_{n-2}\cup\left(\left\langle d_{2n-4} \right\rangle *\bigcup_{i=1}^{n-2}\mathcal{A}_i''\right)
    \end{equation}
    where $\mathcal{A}_i'$ is as before, and $\mathcal{A}_i''$ is obtained by removing $d_{2n-4}$ from every facet of $\mathcal{A}_i'$, for $1 \leq i \leq n-2$.

    When $n =3$, $Z_3 = \Delta_1 \cup \langle d_2 \rangle$. Since both of these are contractible and disjoint, we get $Z_3 \simeq \mathbb{S}^0$. Therefore, $$\Delta_3 \simeq \Sigma^2 Z_3 \simeq \mathbb{S}^2.$$
   
    Observe that, for $n\geq 4$, $\displaystyle{\bigcup_{i=1}^{n-2}}\mathcal A_i''$ is precisely $Z_{n-1}$. 
    Therefore, the expression for $Z_n$, for $n \geq 4$ simplifies to 
    \begin{equation} \label{eqn: Zn and Delta_n together}
        Z_n = \Delta_{n-2} \cup \left( \left\langle d_{2n-4} \right\rangle * Z_{n-1} \right) \mbox{ and } \Delta_{n} \simeq \Sigma^2 Z_n.
    \end{equation}

    \noindent We now examine the link and deletion of $d_{2n-4}$ in $Z_n$. 
    Similar argument used before for $\lk_{X_n}(d_{2n-2})$ and $\del_{X_n}(d_{2n-2})$ implies that $\lk_{Z_n}(d_{2n-4}) = Z_{n-1}$ and $\del_{Z_n}(d_{2n-4}) = \Delta_{n-2}$, respectively. 

   From \Cref{eqn: descriptive form of Zn} and \Cref{eqn: short expression of Xn}, it follows that $Z_{n-1}$ is a subcomplex of $X_{n-2}$.
   % Thus the inclusion of $Z_{n-1}$ is nullhpmotopic in $\st_{\Delta_{n-2}}(a_{2n-4})$. 
    Since the inclusion of $\st_{\Delta_{n-2}}(a_{2n-4})$ in $\Delta_{n-2}$ is null-homotopic, this implies that composition of inclusions 
    \[Z_{n-1} \hookrightarrow X_{n-2} \hookrightarrow \st_{\Delta_{n-2}}(a_{2n-4}) \hookrightarrow \Delta_{n-2}\]
    is null-homotopic.
    However, $\Delta_{n-2} = \del_{Z_n}(d_{2n-4})$. Therefore, $\lk_{Z_n}(d_{2n-4}) \hookrightarrow \del_{Z_n}(d_{2n-4})$ is null-homotopic. Then, by \Cref{lemma: link-deletion}, we have 
    \begin{equation}\label{eqn: Zn recursive homotopy}
        Z_n \simeq \del_{Z_n}(d_{2n-4}) \vee \Sigma \lk_{Z_n}(d_{2n-4}) \simeq \Delta_{n-2} \vee \Sigma Z_{n-1}.
    \end{equation}

    \noindent Therefore, using \Cref{eqn: Zn recursive homotopy} and \Cref{eqn: Zn and Delta_n together}, we get\begin{align*}
        \Delta_n 
        &\simeq \Sigma^2  Z_n \\
        &\simeq \Sigma^2 \left( \Delta_{n-2} \vee \Sigma  Z_{n-1} \right)\\
        &\simeq \Sigma^2 \Delta_{n-2} \vee \Sigma (\Sigma^2 Z_{n-1})\\
        &\simeq \Sigma^2  \Delta_{n-2}  \vee \Sigma  \Delta_{n-1}.
    \end{align*}

    \noindent The above expression along with $\Delta_2 \simeq \mathbb{S}^1$ and $\Delta_3 \simeq \mathbb{S}^2$ implies that for $n\geq 2$, we have $$\Delta_n\simeq\bigvee_{f_{n-1}}\mathbb{S}^{n-1}.$$

    \noindent This completes the proof. 
\end{proof}

\section*{Concluding remarks}
The appearance of Fibonacci numbers in our result is interesting in view of their connection with perfect matchings of ladder graphs. It is well known that the perfect matchings of the $2\times n$ grid graph are in bijective correspondence with domino tilings of a $2\times n$ rectangular board, and the number of all the domino tilings is given by the Fibonacci number $f_{n+1}$. Here, each vertex of the grid graph is viewed as the centre of a square in the board. However, our result reveals a different occurrence of Fibonacci numbers, as they determine the number of spheres appearing in the homotopy decomposition of the perfect matching complex of the $3 \times 2n$ grid graph. This makes the result particularly interesting, as it shows that Fibonacci numbers arise not only in counting perfect matchings, but also in describing the topology of the associated complex. This leads to the following question.

\begin{question}
    Which graph theoretic property leads to the occurrence of the Fibonacci sequence in the homotopy type of the perfect matching complex of the $3 \times 2n$ grid graph? 
\end{question}

%%% ACKNOWLEDGEMENT %%%
\section*{Acknowledgements}

The project started at the Frontiers of Combinatorics: Topology, Algebra and Data‑Driven Insights conference hosted by the Birla Institute of Technology $\&$ Science (BITS), Pilani. We would like to thank BITS Pilani for providing us with wonderful facilities leading to a productive stay. Chandrakar also acknowledges Dr. Anurag Singh for his helpful suggestions during the early stages of this work.

%Chandrakar is partially supported by a grant from the Infosys Foundation. Chau is partially supported in part by the Pacific Institute for the Mathematical Sciences and a grant from the Simons Foundation International [SFI-MPS-T-Institutes-00020822, OY]. Goyal is partially supported by the NFSG grant N5/25/1045.

%%% REFERENCES %%%
\nocite{}
\printbibliography 
\end{document}